\documentclass{amsart}
\usepackage{a4wide,amssymb}

\usepackage[all]{xy}

\usepackage{graphicx}

\makeatletter
\DeclareRobustCommand{\cev}[1]{%
  \mathpalette\do@cev{#1}%
}
\newcommand{\do@cev}[2]{%
  \fix@cev{#1}{+}%
  \reflectbox{$\m@th#1\vec{\reflectbox{$\fix@cev{#1}{-}\m@th#1#2\fix@cev{#1}{+}$}}$}%
  \fix@cev{#1}{-}%
}
\newcommand{\fix@cev}[2]{%
  \ifx#1\displaystyle
    \mkern#23mu
  \else
    \ifx#1\textstyle
      \mkern#23mu
    \else
      \ifx#1\scriptstyle
        \mkern#22mu
      \else
        \mkern#22mu
      \fi
    \fi
  \fi
}

\makeatother

\newcommand{\IR}{\mathbb R}
\newcommand{\IC}{\mathbb C}
\newcommand{\IT}{\mathbb T}
\newcommand{\Aff}{\mathrm{Aff}}
\newcommand{\sign}{\,\mathrm{sign}}

\newcommand{\trig}{\triangle}
\newcommand{\unbal}{K}

\newtheorem{problem}{Problem}
\newtheorem{theorem}{Theorem}
\newtheorem{remark}{Remark}
\newtheorem{lemma}{Lemma}
\newtheorem{corollary}{Corollary}

\theoremstyle{definition}
\newtheorem{definition}{Definition}

\title[$G$-deviations and their applications]{$G$-deviations of polygons and
 their applications\\ in Electric Power Engineering}
\author{Taras Banakh, Olena Hryniv, Vasyl Hudym}
\address{T.Banakh: Ivan Franko National University of Lviv, Ukraine, and Jan Kochanowski University in Kielce, Poland}
\email{t.o.banakh@gmail.com}
\address{O.Hryniv: Ivan Franko National University of Lviv, Ukraine}
\email{ohryniv@gmail.com}
\address{V.Hudym: Faculty of Mechanics and Energetics, Lviv National Agrarian University, Lviv, Ukraine}
\email{vasylhudym@ukr.net}
\subjclass{51M04}
\keywords{Triangle, deviation, group action, Electric Power Engineering}

\begin{document}
\begin{abstract} For any metric space $X$ endowed with the action of a group $G$, and two $n$-gons $\vec x=(x_1,\dots,x_n)\in X^n$ and $\vec y=(y_1,\dots,y_n)\in X^n$ in $X$, we introduce the $G$-deviation $d(G\vec x,\vec y\,)$ of $\vec x$ from $\vec y$ as the distance in $X^n$ from $\vec y$ to the $G$-orbit $G\vec x$ of $\vec x$ in the $n$-th power $X^n$ of $X$. For some groups $G$ of affine transformations of the complex plane, we deduce simple-to-apply formulas for calculating the $G$-deviation between $n$-gons on the complex plane. We apply these formulas for defining new measures of asymmetry of triangles. These new measures can be applied in Electric Power Engineering for evaluating the quality of 3-phase electric power.  One of such measures, namely the affine deviation, is espressible via the unbalance degree, which is a standard characteristic of quality of three-phase electric power.
\end{abstract}
\maketitle

\section*{Introduction}

This paper is motivated by the problem of evaluation and control of the quality of three-phase electric power, which can be represented by a triangle in the Euclidean plane. The vertices of this triangle correspond to the three phases of the electric power, and the lengths of sides of the triangle correspond to the voltage between the phases, and those can be measured by standard volt-meters. So, the triangle is given up to a rotation or even isometry. The quality of electric power can be evaluated by the deviation of the triangle from a regular triangle of a given size. This motivates a mathematical problem of defining a deviation of a polygon in a $G$-space from another polygon. 

This problem can be formalized as follows. Let $X$ be a metric space endowed with an action of some group $G$ whose neutral element is denoted by $1_G$. The action is a function $G\times X\to X$, $(g,x)\mapsto g\cdot x$, satisfying two conditions:
\begin{itemize}
\item $1_X\cdot x=x$ for any $x\in X$;
\item $g\cdot(h\cdot x)=(gh)\cdot x$ for any $g,h\in G$ and $x\in X$.
\end{itemize}

By an {\em $n$-gon} in the metric space $X$ we understand any $n$-tuple $(x_1,\dots,x_n)$ of points of the space $X$. So, $n$-gons are elements of the $n$-th power $X^n$ of the space $X$. The distance between two $n$-gons $\vec x=(x_1,\dots,x_n)$ and $\vec y=(y_1,\dots,y_n)$ in $X$ is calculated by the formula
$$d(\vec x,\vec y\,)=\Big(\sum_{k=1}^n |x_k-y_k|^2\Big)^{\frac12},$$
where $|x_k-y_k|$ is the distance between the points $x_k,y_k$ in the metric space $X$.

The action $G\times X\to X$ induces the coordinatewise action 
$$G\times X^n\to X^n,\quad (g,(x_1,\dots,x_n))\mapsto (gx_1,\dots,gx_n),$$of $G$ on the $n$th power $X^n$ of $X$.

\begin{definition} The {\em $G$-deviation} of an $n$-gon $\vec x=(x_1,\dots,x_n)\in X^n$ from an $n$-gon $\vec y=(y_1,\dots,y_n)\in X^n$ is defined as the real number
$$d(G\vec x,\vec y\,)=\inf\{d(g\vec x,\vec y\,):g\in G\}.$$
\end{definition}
So the $G$-deviation $d(G\vec x,\vec y)$ of $\vec x$ from $\vec y$ is the distance from $\vec y$ to the the $G$-orbit $G\vec x=\{g\vec x:g\in G\}$ of $\vec x$  in the metric space $X^n$. 

In this paper we derive simple formulas for calculating deviations of $n$-gons on the complex plane $\IC$ endowed with the standard Euclidean distance and the action of some groups of affine transformations of $\IC$.
Observe that the subsets
$$\IT=\{z\in\IC:|z|=1\}\quad\mbox{and}\quad\IC^*=\{z\in\IC:z\ne0\}$$ of the complex plane $\IC$ are groups with respect to the operation of multiplication of complex numbers.

A bijective self-map $f$ of the complex plane $\IC$ is called an {\em affine transformation} of $\IC$ if there are complex numbers $a\in\IC^*$ and $b\in\IC$ such that $f(z)=az+b$ for all $z\in\IC$. Geometrically, affine maps of $\IC$ are orientation-preserving similarity transformations of the plane. An affine transformation $az+b$ of $\IC$ is called 
\begin{itemize}
\item {\em linear} if $b=0$;
\item an {\em isometry} if $a\in\IT$;
\item a {\em rotation} if $a\in\IT$ and $b=0$.
\end{itemize}

It is clear that affine transformations of the complex plane form a group with respect to the operation of composition of transformations. Linear transformations, isometries and rotations form  subgroups in the group $\Aff(\IC)$ of affine transformations of $\IC$. Depending on the choice of a subgroup of $\Aff(\IC)$, we distinguish four deviations of an $n$-gon $\vec x\in\IC^n$ from an $n$-gon $\vec y\in\IC^n$:
\begin{itemize}
\item the {\em rotational deviation} $d(\IT\vec x,\vec y\,)=\inf\{d(a\vec x,\vec y\,):a\in\IT\}$;
\item the {\em linear deviation} $d(\IC^*\vec x,\vec y\,)=\inf\{d(a\vec x,\vec y\,):a\in\IC^*\}$;
\item the {\em isometric deviation} $d(\IT\vec x+\IC,\vec y\,)=\inf\{d(a\vec x+b,\vec y\,):a\in\IT,\;b\in\IC\}$;
\item the {\em affine deviation} $d(\IC^*\vec x+\IC,\vec y\,)=\inf\{d(a\vec x+b,\vec y\,):a\in\IC^*,\;b\in\IC\}$.
\end{itemize}
In Sections~\ref{s2}--\ref{s4} we shall deduce simple formulas for calculating these four deviations for $n$-gons in the complex plane. 

In Section~\ref{s5} we apply these formulas for calculating deviations between triangles in the complex plane. 

By a {\em triangle} in a metric space $X$ we understand any triple $\vec x=(x_1,x_2,x_2)$ of points of the space $X$. 
For two points $x,y\in X$ we denote by $|x-y|$ the distance between $x$ and $y$ in the metric space $X$.

A triangle $\vec x=(x_1,x_2,x_3)$ in a metric space $X$ is defined to be
\begin{itemize}
\item {\em equilateral} if $|x_1-x_2|=|x_2-x_3|=|x_3-x_1|$;
\item {\em singular} if $x_1=x_2=x_3$;
\item {\em regular} if it is equilateral and not singular;
\item {\em linear} if there exist pairwise distinct numbers $i,j,k\in\{1,2,3\}$ such that $|x_i-x_k|=|x_i-x_j|+|x_j-x_k|$.
\end{itemize}
Observe that a triangle is linear if and only if it is isometric to a triangle in the real line. A triangle is singular if and only if it is both linear and equilateral.

The main result of Section~\ref{s5} is Theorem~\ref{t:main} providing formulas for calculating the affine and isometric deviations of a plane triangle (given by lengths of its sides) from the regular triangle
$$\trig=\tfrac1{\sqrt{3}}(1,e^{i\frac{2\pi}3},e^{-i\frac{2\pi}3})$$with sides $1$. The triangle $\trig$ is oriented counter-clockwise  since the move from its first vertex $1$ to the second one $e^{i\frac{2\pi}3}$ goes counter-clockwise. For a triangle $\vec z$ in the complex plane its {\em orientation} is defined as
$$\sign(\vec z\,)=\begin{cases}
-1&\mbox{if $\vec z$ is singular};\\
0&\mbox{if $\vec z$ is linear and not singular};\\
1&\mbox{if $\vec z$ is non-linear and is oriented counter-clockwise};\\
-1&\mbox{if $\vec z$ is non-linear and is oriented clockwise}.
\end{cases}
$$
According to this definition, any singular triangle has opposite orientation to the regular triangle $\triangle$. The regular triangle $\trig$ has orientation $\sign(\trig)=1$.

In Theorem~\ref{t:main} we prove that for a triangle $\vec z=(z_1,z_2,z_3)$ on the complex plane, 
\begin{itemize}
\item the isometric deviation $d(\IT\vec z+\IC,\trig)=\sqrt{1+u^2-\sqrt{2}u\sqrt{1+\sign(\vec z\,)\sqrt{3-6q}}}$,
\item the affine deviation $d(\IC^*\vec z+\IC,\trig)=\sqrt{\frac12\big(1-\sign(\vec z\,)\sqrt{3-6q}\,\big)}$,
\end{itemize}
where $$u=\sqrt{\frac{|z_1-z_2|^2+|z_2-z_3|^2+|z_3-z_1|^2}3}\quad\mbox{and}\quad q=\frac{|z_1-z_2|^4+|z_2-z_3|^4+|z_3-z_1|^4}{(|z_1-z_2|^2+|z_2-z_3|^2+|z_3-z_1|^2)^2}.$$
If $z_1=z_2=z_3$, then we put $q=\frac12$.

The coefficient $q$ is called the {\em quadrofactor} of the triangle $\vec z$. It takes its minimal value $\frac 13$ if and only if the triangle is regular, and $q$ takes its maximal value $\frac12$ if and only if the triangle is linear, see Theorem~\ref{t:quadro}. 

The number $\sqrt{3-6q}$ appearing in the above formulas has a nice geometric meaning: it is equal to $\frac{4A}{\sqrt{3}u^2}$ where $A$ is the area of the triangle $\vec z$. The fraction $\frac{4A}{\sqrt{3}u^2}=\sqrt{3-6q}$ is called the {\em normalized area} of the triangle $\vec z$. It can be considered as a measure of regularity of the triangle: the normalized are takes its
\begin{itemize}
\item maximal value $1$ if and only if the triangle is regular, 
\item minimal value $0$ if and only if the triangle is linear,
\end{itemize}
see Corollary~\ref{c:quadro}.

The formula for the affine deviation and the properties of the quadrofactor imply that for a triangle $\vec z$ on the complex plane its affine deviation $d(\IC^*\!\vec z{+}\IC,\trig)$ from the regular triangle $\trig$ takes
\begin{itemize}
\item its minimal value $0$ if and only if the triangle $\vec z$ is regular and oriented counter-clockwise;
\item its maximal value $1$ if and only if the triangle $\vec z$ is equilateral and oriented clockwise;
\item the intermediate value $\frac1{\sqrt{2}}$ if and only if the triangle $\vec z$ is linear but nor singular;
\item its value in the interval $(0,\frac1{\sqrt{2}})$ if and only if the triangle $\vec z$ is not linear, not regular, and is oriented counter-clockwise;
\item its value in the interval $(\frac1{\sqrt{2}},1)$ if and only if the triangle $\vec z$ is not linear, not regular, and is oriented clockwise.
\end{itemize}
  
We propose to use the isometric and affine deviations as characteristics of the quality of 3-phase electric power. In Section~\ref{s:K} we compare these characteristics to a standard measure of quality of electric power, which is based on positive and negative components of a triangle. Following the classical approach of Fortescue \cite{Fort}, for a triangle $\vec z=(z_1,z_2,z_3)$ in the complex plane, we define three symmetric components of $\vec z$:
\begin{itemize}
\item the {\em positive component} $\frac13(z_1+z_2e^{-\frac{2\pi}3i}+z_3e^{\frac{2\pi}3i})$,
\item the {\em negative component} $\frac13(z_1+z_2e^{\frac{2\pi}3i}+z_3e^{-\frac{2\pi}3i})$, and
\item the {\em zero component} $\frac13(z_1+z_2+z_3)$.
\end{itemize}
The ratio of negative to positive components
$$\unbal=\frac{|z_1+z_2e^{\frac{2\pi}3i}+z_3e^{-\frac{2\pi}3i}|}{|z_1+z_2e^{-\frac{2\pi}3i}+z_3e^{\frac{2\pi}3i}|}$$
is called the {\em unbalance factor} of the triangle $\vec z$ and is widely used in Electric Power Engineering as a number characterizing a deviation of the triangle $\vec z$ from being regular and oriented counter-clockwise. The unbalance factor $\unbal$ is equal to zero if and only if the triangle $\vec z$ is regular and is oriented counter-clockwise. In particular, the unbalance factor $\unbal$ of the model regular triangle $\trig$ is zero.

In Theorem~\ref{t:K} and Corollary~\ref{c:K} we prove that any triangle $\vec z=(z_1,z_2,z_3)\in\IC$ has unbalance factor    
$$\unbal=\sqrt{\frac{1-\sign(\vec z\,)\sqrt{3-6q}}{1+\sign(\vec z\,)\sqrt{3-6q}}}=\sqrt{\frac{d(\IC^*\vec z+\IC,\trig)^2}{1-d(\IC^*\vec z+\IC,\trig)^2}},$$ which implies that $\unbal$ is expressible via the affine deviation $$d(\IC^*\vec z+\IC,\trig)=\sqrt{\tfrac12(1-\sign(\vec z\,)\sqrt{3-6q})}=\sqrt{\frac{\unbal^2}{1+\unbal^2}}$$ and vice versa.

On the other hand, the isometric deviation $d(\IT\vec z+\IC,\trig)$ is a new characteristic that evaluates a deviation of a triangle $\vec z$ from $\trig$ both by the form and size.

\section{The Hilbert space structure of $\IC^n$}\label{s1}

For a complex number $z=x+iy$  we denote by $|z|=\sqrt{x^2+y^2}$ the absolute value of $z$. Also $\Re(z)$ and $\Im(z)$ will denote the real and imaginary parts of $z=x+iy$, which are equal to $x$ and $y$, respectively. By $\bar z:=x-iy$ we denote the conjugate to $z$. We shall often exploit the equality $|z|^2=z\bar z$ holding for any complex number $z$.
\smallskip

For two vectors $\vec x=(x_1,\dots,x_n)$ and $\vec y=(y_1,\dots,y_n)$ in the linear space $\IC^n$ by
$$\langle \vec x|\vec y\,\rangle:=\sum_{k=1}^n{x_k\overline{y_k}}$$we denote the scalar product of the vectors $\vec x$ and $\vec y$. The real number $$\|\vec x\|:=\sqrt{\langle \vec x |\vec x\,\rangle}=\Big(\sum_{k=1}^n{|x_k|^2}\Big)^{\frac12}$$is called the {\em norm} of the vector $\vec x$ in $\IC^n$.

It is easy to see that for two $n$-gons $\vec x=(x_1,\dots,x_n)$ and $\vec y=(y_1,\dots,y_n)$ in $\IC$ we have
$$d(\vec x,\vec y\,)=\Big(\sum_{k=1}^n{|x_k-y_k|^2}\Big)^\frac12=\|\vec x-\vec y\,\|.$$
So, the distance between $n$-gons is equal to the distance between the corresponding vectors in the Hilbert space $\IC^n$.

For an $n$-gon  $\vec x=(x_1,\dots,x_n)$ in $\IC$ and a complex number $z$ let
$$\vec x+z=(x_1+z,\dots,x_n+z)$$be the shifted $n$-gon $\vec x$ in the direction of the vector $z$.

\section{The rotational deviation}

\begin{theorem}\label{t1}  The rotational deviation $d(\IT\vec x,\vec y\,)$ of any vector $\vec x=(x_1,\dots,x_n)\in\IC^n$ from a vector $\vec y=(y_1,\dots,y_n)\in\IC^n$ can be  calculated by the formula
$$d(\IT\vec x,\vec y\,)=\sqrt{\|\vec x\|^2+\|\vec y\,\|^2-2\cdot|\langle \vec x\,|\vec y\,\rangle|}=d(a\vec x,\vec y\,),$$
where $a=\frac{\overline{\langle \vec x|\vec y\,\rangle}}{|\langle \vec x|\vec y\,\rangle|}$ if $\langle \vec x|\vec y\,\rangle\ne 0$, and $a$ is any element of $\IT$ if $\langle \vec x|\vec y\,\rangle=0$.\end{theorem}

\begin{proof} 
Observe that for any $a\in\IT$ we have
$$
\begin{aligned}
d(a\vec x,\vec y\,)^2=&\sum_{k=1}^n|ax_k-y_k|^2=\sum_{k=1}^n(ax_k-y_k)(\bar a\bar x_k-\bar y_k)=\\
&\sum_{k=1}^n\big(a\bar a x_k\bar x_k+y_k\bar y_k-(ax_k\bar y_k+\bar a\bar x_ky_k)\big)=\sum_{k=1}^n\big(|x_k|^2+|y_k|^2-2\Re(ax_k\bar y_k)\big)=\\
&\big(\sum_{k=1}^n|x_k|^2+|y_k|^2\big)-2\Re(a\sum_{k=1}^nx_k\bar y_k)=\|\vec x\|^2+\|\vec y\,\|^2-2\Re(a\langle \vec x|\vec y\,\rangle).
\end{aligned}
$$
Now we see that
$$
\begin{aligned}
d(\IT\vec x,\vec y\,)=&\min\{d(a\vec x,\vec y\,):a\in\IT\}=\min\big\{\|\vec x|^2+|\vec y\,\|^2-2\Re(a\langle \vec x|\vec y\,\rangle):a\in\IT\}=\\
&\|\vec x\|^2+\|\vec y\,\|^2-2|\langle \vec x|\vec y\,\rangle|=d(a\vec x,\vec y\,),
\end{aligned}
$$
where $a=\frac{\overline{\langle \vec x|\vec y\,\rangle}}{|\langle \vec x|\vec y\,\rangle|}$ if $\langle \vec x|\vec y\,\rangle\ne 0$, and $a$ is any element of $\IT$ if $\langle \vec x,\vec y\,\rangle=0$.\end{proof}

\section{The isometric deviation}\label{s2}

\begin{theorem}\label{t2}  The isometric deviation $d(\IT\vec x+\IC,\trig)$ of any vector $\vec x=(x_1,\dots,x_n)\in\IC^n$ from a vector $\vec y=(y_1,\dots,y_n)\in\IC^n$  can be  calculated by the formula
$$d(\IT\vec x+\IC,\vec y\,)=\sqrt{\|\vec x-x_0\|^2+\|\vec y-y_0\|^2-2\cdot|\langle \vec x-x_0|\vec y-y_0\rangle|}=d(a(\vec x-x_0),\vec y-y_0),$$
where $$x_0=\frac1n\sum_{k=1}^nx_k,\;\;y_0=\frac1n\sum_{k=1}^ny_k,$$ and $$a=\frac{\overline{\langle \vec x-x_0|\vec y-y_0\rangle}}{|\langle \vec x-x_0|\vec y-y_0\rangle|}$$ if $\langle \vec x-x_0|\vec y-y_0\rangle\ne 0$, and $a$ is any element of $\IT$ if $\langle \vec x-x_0|\vec y-y_0\rangle= 0$.
\end{theorem}

\begin{proof}
Observe that 
$$d(\IT\vec x+\IC,\vec y\,)^2=\min\{d(a(\vec x-x_0)+b,\vec y-y_0)^2:a\in\IT,\;b\in\IC\}.$$
For any $a\in\IT$ and $b\in\IC$ we have
$$
\begin{aligned}
&d(a(\vec x-x_0)+b,\vec y-y_0)^2=\sum_{k=1}^n(a(x_k-x_0)+b-(y_k-y_0))(\bar a(\bar x_k-\bar x_0)+\bar b-(\bar y_k-\bar y_0))=\\
&=\sum_{k=1}^n\big(|a(x_k{-}x_0)-(y_k{-}y_0)|^2+|b^2|+b(\bar a(\bar x_k{-}\bar x_0)-(\bar y_k{-}y_0))+(a(x_k{-}x_0)-(y_k{-}y_0))\bar b\big)=\\
&=\|a(\vec x{-}x_0)-(\vec y{-}y_0)\|^2+|b^2|+b\bar a\sum_{i=1}^n(\bar x_k{-}\bar x_0)-b\sum_{k=1}^n(\bar y_k{-}y_0)+\bar ba\sum_{i=1}^n({x_k{-}x_0})-\bar b\sum_{k=1}^n(y_k{-}y_0)=\\
&=d(a(\vec x-x_0),\vec y-y_0)^2+|b|^2+b\bar a\cdot 0-b\cdot 0+\bar ba\cdot 0-\bar b\cdot0=d(a(\vec x-x_0),\vec y-y_0)^2+|b|^2.
\end{aligned}
$$

Applying Theorem~\ref{t1}, we obtain
$$
\begin{aligned}
d(\IT\vec x+\IC,\vec y\,)^2=&\min\{d(a(\vec x-x_0)+b,\vec y-y_0)^2:a\in\IT,\;b\in\IC\}=\\
&\min\{d(a(\vec x-x_0),\vec y-y_0)^2+|b|^2:a\in\IT,\;b\in\IC\}=\\
&\min\{d(a(\vec x-x_0),\vec y-y_0)^2:a\in\IT\}=d(\IT(\vec x-x_0),\vec y-y_0)^2=\\
&\|\vec x-x_0\|^2+\|\vec y-y_0\|^2-2|\langle\vec x-x_0|\vec y-y_0\rangle|=d(a(\vec x-x_0),\vec y-y_0)^2,
\end{aligned}
$$where
$$a=\frac{\overline{\langle \vec x-x_0|\vec y-y_x\rangle}}{|\langle \vec x-x_0|\vec y-y_0\rangle|}$$ if $\langle \vec x-x_0|\vec y-y_0\rangle\ne 0$, and $a$ is any element of $\IT$ if $\langle \vec x-x_0|\vec y-y_0\rangle= 0$.
\end{proof}

\section{The linear deviation}

\begin{theorem}\label{t3}  The linear deviation $d(\IC^*\!\vec x,\vec y\,)$ of any non-zero vector $\vec x=(x_1,\dots,x_n)\in\IC^n$ from a vector $\vec y=(y_1,\dots,y_n)\in\IC^n$ can be  calculated by the formula
$$d(\IC^*\vec x,\vec y\,)=\sqrt{\|\vec y\,\|^2-\frac{|\langle\vec x|\vec y\,\rangle|^2}{\|\vec x\|^2}}=d(a\vec x,\vec y\,),$$
where $$a=\frac{\overline{\langle \vec x|\vec y\,\rangle}}{\|\vec x\|^2}.$$
If $\vec x=\vec 0$, then $d(\IC^*\!\vec x,\vec y\,)=\|\vec y\|$.
\end{theorem}

\begin{proof} If $\vec x=\vec 0$, then 
$$d(\IC^*\!\vec x,\vec y\,)=\min\{d(a\vec 0,\vec y\,):a\in\IC^*\}=d(\vec 0,\vec y\,)=\|\vec y\,\|.$$
So, we assume that $\vec x\ne\vec 0$.

Observe that
$$d(\IC^*\!\vec x,\vec y\,)=\min\{d(tr\vec x,\vec y\,):t\in\IT,\;r\in\IR_+\},
$$where $\IR_+$ stands for the set of positive real numbers.

By Theorem~\ref{t1}, for any $r\in\IR_+$
$$
\min\{d(tr\vec x,\vec y\,)^2:t\in\IT\}=r^2\|\vec x\|^2+\|\vec y\,\|^2-2r|\langle \vec x|\vec y\,\rangle|=
\Big(r\|\vec x\|-\frac{|\langle\vec x|\vec y\,\rangle|}{\|\vec x\|}\Big)^2+\|\vec y\,\|^2-\frac{|\langle\vec x|\vec y\,\rangle|^2}{\|x\|^2}.
$$
Now we see that
$$\begin{aligned}
&d(\IC^*\!\vec x,\vec y\,)^2=\|\vec y\,\|^2-\frac{|\langle\vec x|\vec y\,\rangle|^2}{\|x\|^2}=d(tr\vec x,\vec y\,)^2=d(a\vec x,\vec y\,)^2,
\end{aligned}
$$
where  $r=\frac{|\langle\vec x|\vec y\,\rangle|}{\|\vec x\|^2}$, 
$t=\frac{\overline{\langle r\vec x|\vec y\,\rangle}}{|\langle r\vec x|\vec y\,\rangle|}=\frac{\overline{\langle \vec x|\vec y\,\rangle}}{|\langle \vec x|\vec y\,\rangle|}
$, and $a=\frac{\overline{\langle\vec x|\vec y\,\rangle}}{\|\vec x\|^2}$.
\end{proof}

\section{The affine deviation}\label{s4}

\begin{theorem}\label{t4}  The affine deviation $d(\IC^*\!\vec x+\IC,\vec y\,)$ of any non-constant vector $\vec x=(x_1,\dots,x_n)\in\IC^n$ from a vector $\vec y=(y_1,\dots,y_n)\in\IC^n$ can be  calculated by the formula
$$d(\IC^*\!\vec x+\IC,\vec y\,)=\sqrt{\|\vec y-y_0\|^2-\frac{|\langle\vec x-x_0|\vec y-y_0\rangle|^2}{\|\vec x-x_0\|^2}}=d(a(\vec x-x_0),\vec y-y_0),$$
where $$a=\frac{\overline{\langle \vec x-x_0|\vec y-y_0\rangle}}{\|\vec x-x_0\|^2}.$$
If $\vec x\in\IC^n$ is a constant vector, then $d(\IC^*\vec x+\IC,\vec y\,)=\|\vec y-y_0\|$ for any vector $\vec y\in\IC^n$.
\end{theorem}

\begin{proof}
Observe that 
$$d(\IC^*\vec x+\IC,\vec y\,)^2=\min\{d(a(\vec x-x_0)+b,\vec y-y_0)^2:a\in\IC^*,\;b\in\IC\}.$$
For any $a\in\IC^*$ and $b\in\IC$ we have
$$
\begin{aligned}
&d(a(\vec x-x_0)+b,\vec y-y_0)^2=\sum_{k=1}^n(a(x_k-x_0)+b-(y_k-y_0))(\bar a(\bar x_k-\bar x_0)+\bar b-(\bar y_k-\bar y_0))=\\
&=\sum_{k=1}^n\big(|a(x_k-x_0)-(y_k-y_0)|^2+|b^2|+b(\bar a(\bar x_k-\bar x_0)-(\bar y_k-\bar y_0))+(a(x_k-x_0)-(y_k-y_0))\bar b\big)=\\
&=\|a(\vec x{-}x_0)-(\vec y{-}y_0)\|^2+|b^2|+b\bar a\sum_{i=1}^n(\bar x_k{-}\bar x_0)-b\sum_{k=1}^n(\bar y_k{-}\bar y_0)+\bar ba\sum_{i=1}^n(x_k{-}x_0)-\bar b\sum_{k=1}^n(y_k{-}y_0)=\\
&=d(a(\vec x-x_0),\vec y-y_0)^2+|b|^2+b\bar a\cdot 0-b\cdot 0+\bar ba\cdot 0-\bar b\cdot0=d(a(\vec x-x_0),\vec y-y_0)^2+|b|^2.
\end{aligned}
$$
Then
$$
\begin{aligned}
d(\IC^*\!\vec x+\IC,\vec y\,)^2=&\min\{d(a(\vec x-x_0)+b,\vec y-y_0)^2:a\in\IC^*,\;b\in\IC\}=\\
&\min\{d(a(\vec x-x_0),\vec y-y_0)^2+|b|^2:a\in\IC^*,\;b\in\IC\}=\\
&\min\{d(a(\vec x-x_0),\vec y-y_0)^2:a\in\IC^*\}=d(\IC^*(\vec x-x_0),\vec y-y_0)^2.
\end{aligned}
$$
If $\vec x$ is a contant vector, then $\vec x-x_0=\vec 0$ and 
$$d_1(\IC^*\!\vec x+\IC,\vec y\,)=d(\IC^*(\vec x-x_0),\vec y-y_0)=d(\IC^*{\cdot}\vec 0,\vec y-y_0)=\|\vec y-y_0\|.$$
If the vector $\vec x$ is not constant, then $\vec x-x_0\ne\vec 0$ and by  Theorem~\ref{t3}, we obtain
$$
\begin{aligned}
d(\IC^*\!\vec x+\IC,\vec y\,)^2=d(\IC^*(\vec x-x_0),\vec y-y_0)=\|\vec y-y_0\|^2-\frac{|\langle\vec x-x_0|\vec y-y_0\rangle|^2}{\|\vec x-x_0\|^2}=d(a(\vec x-x_0),\vec y-y_0),
\end{aligned}
$$where
$a=\frac{\overline{\langle \vec x-x_0|\vec y-y_0\rangle}}{\|\vec x-x_0\|^2}$.
\end{proof}

\section{The quadrofactor and normalized area of a triangle}

By a {\em plane triangle} we understand a non-singular triangle in the complex plane.

\begin{definition} For a plane triangle $\vec z$ with sides $a,b,c$ in the plane the {\em quadrofactor} of  $\vec z$ is the number
$$q=\frac{a^4+b^4+c^4}{(a^2+b^2+c^2)^2}$$
and the {\em normalized area} of $\vec z$ is the number
$$\frac{4\sqrt{3}A}{a^2+b^2+c^2}$$where $A$ is the area of the triangle $\vec z$.
\end{definition}

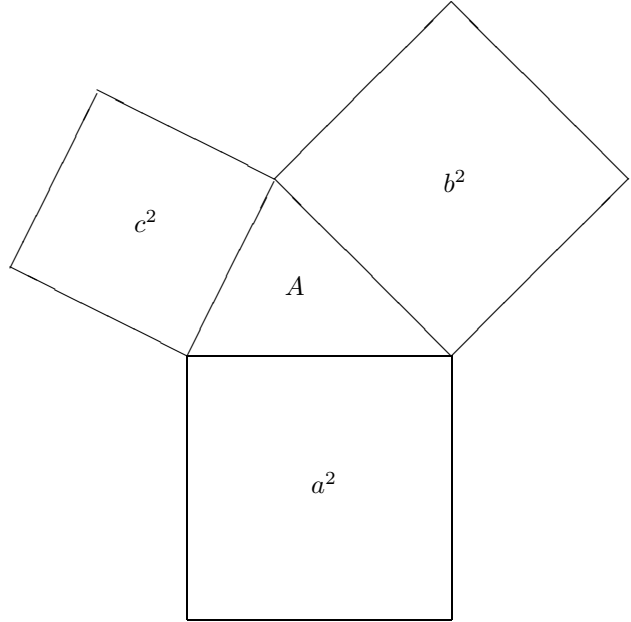
\begin{picture}(300,250)(-50,-10)

\put(100,0){\line(1,0){100}}
\put(100,100){\line(1,0){100}}
\put(100,0){\line(0,1){100}}
\put(200,0){\line(0,1){100}}

\put(147,48){$a^2$}
\put(197,162){$b^2$}
\put(80,147){$c^2$}
\put(137,123){$A$}

\put(100,100){\line(1,2){33}}
\put(200,100){\line(-1,1){67}}

\put(200,100){\line(1,1){67}}
\put(133,167){\line(1,1){67}}
\put(267,167){\line(-1,1){67}}

\put(100,100){\line(-2,1){67}}
\put(133,167){\line(-2,1){67}}
\put(33,133){\line(1,2){33}}

\end{picture}

The normalized area can be expressed via the quadrofactor as follows.

\begin{theorem}\label{t:Heron} For a plane triangle with sides $a,b,c$ its normalized area equals 
$$\frac{4\sqrt{3}A}{a^2+b^2+c^2}=\sqrt{3-6q},$$
where $A$ is the area of the triangle and $q$ is its quadrofactor.
\end{theorem}

\begin{proof} By Heron's formula \cite{Heron},
$$
A=\tfrac14\sqrt{2a^2b^2+2a^2c^2+2b^2c^3-a^4-b^4-c^4}=\tfrac14\sqrt{(a^2+b^2+c^2)^2-2(a^4+b^4+c^4)}
$$
and hence
$$\frac{4\sqrt{3}A}{(a^2+b^2+c^2)}=\frac{\sqrt{3}\sqrt{(a^2+b^2+c^2)^2-2(a^4+b^4+c^4)}}{(a^2+b^2+c^2)}=\sqrt{3-6q}
$$
\end{proof}

The following theorem describes three basic properties of the quadrofactor.

\begin{theorem}\label{t:quadro} Let $q$ be the quadrofactor of a plane triangle. Then
\begin{enumerate}
\item $\frac13\le q\le\frac12$;
\item $q=\frac13$ if and only if the triangle  is regular;
\item $q=\frac12$ if and only if the triangle  is linear.
\end{enumerate}
\end{theorem}

\begin{proof} Let $a,b,c$ be the lengths of the sides of the triangle. 
By the Power Mean Inequality \cite[3.1.1]{Hand}, we have
$$\sqrt[2]{\frac{a^2+b^2+c^2}3}\le \sqrt[4]{\frac{a^4+b^4+c^4}3}$$and hence
$$(a^2+b^2+c^2)^2\le 3(a^4+b^4+c^4),$$which implies $q\ge\frac13$. By \cite[3.1.1]{Hand} the equality is attained if and only if $a=b=c$ if and only if the triangle is regular.

Theorem~\ref{t:Heron} implies that $3-6q\ge 0$ and hence $q\le \frac12$. Moreover, $q=\frac12$ if and only if the triangle has zero area if only if the triangle is linear.
\end{proof}

\begin{corollary}\label{c:quadro} Let $q$ be the quadrofactor of a plane triangle and $\sqrt{3-6q}$ be its normalized area. Then
\begin{enumerate}
\item $\sqrt{3-6q}\le 1$;
\item $\sqrt{3-6q}=1$ if and only if the triangle is regular;
\item $\sqrt{3-6q}=0$ if and only if the triangle  is linear.
\end{enumerate}
\end{corollary}

\section{Deviations from a regular triangle}\label{s5}

In this section we deduce formulas for calculating four deviations of a given triangle on the complex plane from the regular triangle
$$\trig=\tfrac1{\sqrt{3}}(1,e^{i\frac{2\pi}3},e^{-i\frac{2\pi}3})$$
with sides of length $1$. This regular triangle will model the ideal form of a 3-phase electric power.

Theorems~\ref{t1}--\ref{t4} imply the following formulas for deviations of a triangle from  the regular triangle $\trig$.

\begin{theorem}\label{t5} A triangle $\vec z=(z_1,z_2,z_3)\in\IC^3$ with center $ z_0=\tfrac13(z_1+z_2+z_3)$ has:
\begin{enumerate}
\item the rotation derivation $d(\IT\vec z,\trig)=\sqrt{\|\vec z\|^2+1-2\cdot|\langle \vec z|\trig\rangle|}=d\big(\frac{\overline{\langle \vec z\,|\triangle\rangle}}{|\langle \vec z\,|\triangle\rangle|}\vec z,\trig\big)$,
\item the linear deviation $d(\IC^*\vec z,\trig)=\sqrt{1-\frac{|\langle \vec z|\trig\rangle|}{\|\vec z\|^2}}=d\big(\frac{\overline{\langle \vec z\,|\triangle\rangle}}{\|\vec z\,\|^2}\vec z,\trig\big)$,
\item the isometric deviation $d(\IT\vec z{+}\IC,\trig){=}\sqrt{\|\vec z{-}z_0\|^2{+}1{-}2{\cdot}|\langle \vec z{-}z_0|\trig\rangle|}{=}d\big(\frac{\overline{\langle \vec z-z_0|\triangle\rangle}}{|\langle \vec z{-}z_0|\triangle\rangle|}(\vec z{-}z_0),\trig\big)$,
\item the affine deviation $d(\IC^*\vec z+\IC,\trig)=\sqrt{1-\frac{|\langle \vec z-z_0|\trig\rangle|^2}{\|\vec z-z_0\|^2}}=d\big(\frac{\overline{\langle \vec z-z_0\,|\triangle\rangle}}{\|\vec z-z_0\,\|^2}(\vec z-z_0),\trig\big)$.
\end{enumerate}
\end{theorem}

Next, we show that the isometric and affine deviations $d(\IT\vec z+\IC,\trig)$ and $d(\IC^*\vec z+\IC)$ can be expressed via the lengths of the sides of the triangle $\vec z$. 

\begin{lemma}\label{l1} For any triangle $\vec z=(z_1,z_2,z_3)\in\IC^2$ on the complex plane 
and its center $z_0=\frac13(z_0,z_1,z_2)$ we have
$$\|\vec z-z_0\|^2=\tfrac13(|z_1-z_2|^2+|z_2-z_3|^2+|z_3-z_1|^2).$$
\end{lemma}

\begin{proof} Let $x_k=z_k-z_0$ for $k\in\{1,2,3\}$ and observe that $x_1+x_2+x_3=0$.
Observe also that $x_k-x_j=z_k-z_j$ for any $k,j\in\{1,2,3\}$. Then 
$$
\begin{aligned}
&|z_1-z_2|^2+|z_2-z_3|^2+|z_3-z_1|^2=|x_1-x_2|^2+|x_2-x_3|^2+|x_3-x_1|^2=\\
&=(x_1-x_2)(\bar x_1-\bar x_2)+(x_2-x_3)(\bar x_2-\bar x_3)+(x_3-x_1)(\bar x_3-\bar x_1)=\\
&=2(x_1\bar x_1+x_2\bar x_2+x_3\bar x_3)-x_1\bar x_2-x_2\bar x_1-x_2\bar x_3-x_3\bar x_2-x_3\bar x_1-x_1\bar x_3=\\
&=2(|x_1|^2+|x_2|^2+|x_3|^2)-x_1\bar x_2+(x_1+x_3)\bar x_1-x_2\bar x_3+(x_1+x_2)\bar x_2-x_3\bar x_1+(x_2+x_3)\bar x_3=\\
&=3(|x_1|^2+|x_2|^2+|x_3|^2)-x_1\bar x_2+x_3\bar x_1-x_2\bar x_3+x_1\bar x_2-x_3\bar x_1+x_2\bar x_3=\\
&=3\cdot (|z_1-z_0|^2+|z_2-z_0|^2+|z_3-z_0|^2)+0=3\cdot\|\vec z-z_0\|^2.
\end{aligned}
$$
\end{proof}

\begin{lemma}\label{l2} Let $\vec z=(z_1,z_2,z_3)$ where $z_1=0$, $z_2=x_2\in\IR$ and $z_3=x_3+iy_3$. Then for the center $z_0=\frac13(z_1+z_2+z_3)$ of the triangle $\vec z$ we have
$$\langle \vec z-z_0|\trig\rangle=-\frac{x_2+x_3+\sqrt 3 y_3}{2\sqrt{3}}-\frac{\sqrt{3}x_2-\sqrt 3x_3+y_3}{2\sqrt{3}}\cdot i$$
and
$$|\langle \vec z-z_0|\trig\rangle|^2=\tfrac13(x_2^2+x_3^2+y_3^2-x_2x_3+\sqrt{3}x_2y_3).$$
\end{lemma}

\begin{proof} It follows that $$z_0=\tfrac{x_2+x_3}{3}+\tfrac{y_3}{3}i$$and hence
$$\vec z-z_0=(-\tfrac{x_2+x_3}{3}-i\tfrac{y_3}{3}, \tfrac{2x_2-x_3}{3}-i\tfrac{y_3}{3}, \tfrac{2x_3-x_2}{3}+i\tfrac{2y_3}{3}).$$
Since 
$$\trig=\tfrac1{\sqrt{3}}(1,e^{i\tfrac{2\pi}3},e^{-i\tfrac{2\pi}3})=\tfrac1{\sqrt{3}}(1, -\tfrac{1}{2}+i\tfrac{\sqrt{3}}{2}, -\tfrac{1}{2}-i\tfrac{\sqrt{3}}{2}),$$
we have
$$	
	\begin{aligned}
	&\langle \vec z-z_0|\trig\rangle=\tfrac1{\sqrt{3}}\Big((-\tfrac{x_2+x_3}{3}-i\tfrac{y_3}{3})\cdot 1+ (\tfrac{2x_2-x_3}{3}-i\tfrac{y_3}{3})(-\tfrac{1}{2}-i\tfrac{\sqrt{3}}{2})+( \tfrac{2x_3-x_2}{3}+i\tfrac{2y_3}{3})(-\tfrac{1}{2}+i\tfrac{\sqrt{3}}{2})\Big)\\  
	&=\tfrac{1}{3\sqrt{3}}\left(-x_2-x_3-i y_3+ (2x_2-x_3-i y_3)(-\tfrac{1}{2}-i\tfrac{\sqrt{3}}{2})+( 2x_3-x_2+2iy_3)(-\tfrac{1}{2}+i\tfrac{\sqrt{3}}{2}) \right)\\
	&= \tfrac{1}{3\sqrt{3}}\left({-}x_2{-}x_3{-}i y_3- x_2{+}\tfrac{x_3}{2}{+}\tfrac{y_3}{2}i-\sqrt{3}x_2i {+}\tfrac{\sqrt 3}{2}x_3 i{-}\tfrac{\sqrt 3}{2}y_3- x_3{+}\tfrac{x_2}{2}{-}y_3 i+\sqrt{3}x_3i{-}\tfrac{\sqrt 3}{2}x_2 i {-} \sqrt{3}y_3\right) \\
	&= \tfrac{1}{3\sqrt{3}}\left(-\tfrac{3}{2}x_2-\tfrac{3}{2}x_3 -  \tfrac{3\sqrt 3}{2}y_3   +i(\tfrac{-3\sqrt{3}}{2}x_2 +\tfrac{3\sqrt 3}{2}x_3-\tfrac{3}{2} y_3) \right)\\
	&=	\tfrac{1}{2\sqrt{3}}\left(-x_2-x_3 - \sqrt 3 y_3   +i(-\sqrt{3}x_2 +\sqrt 3x_3- y_3) \right).
		\end{aligned}
	$$
	Then
	$$
	\begin{aligned}
	&|\langle \vec z-z_0|\trig\rangle|^2=\tfrac{1}{12}\big( ( -x_2-x_3 -  \sqrt 3 y_3)^2 + (- \sqrt{3}x_2 +\sqrt 3x_3- y_3)^2  \big)\\
	&= \tfrac{1}{12}\big(x_2^2{+}x_3^2 {+} 3 y_3^2 {+}2x_2x_3{+}2\sqrt{3}x_2y_3{+}2\sqrt{3}x_3y_3 + 3x_2^2 {+}3x_3^2{+} y_3^2{-}6x_2x_3{+}2\sqrt{3}x_2y_3{-}2\sqrt{3}x_3y_3  \big)\\
	&= \tfrac{1}{12}\big(4x_2^2+4x_3^2 + 4 y_3^2 -4x_2x_3+4\sqrt{3}x_2y_3 \big).
	 \end{aligned}
	 $$
\end{proof}

\begin{lemma}\label{l3}
For any triangle $\vec z=(z_1,z_2,z_3)\in\IC^3$ and its center $z_0=\frac13(z_1+z_2+z_3)$ we have 
$$
|\langle \vec z-z_0|\trig\rangle|=\frac{u}{\sqrt{2}}\sqrt{1+\sign(\vec z\,)\sqrt{3-6q}}
$$
where $$u=\sqrt{\frac{|z_1-z_2|^2+|z_2-z_3|^2+|z_3-z_1|^2}3}\quad\mbox{and}\quad q=\frac{|z_1-z_2|^4+|z_2-z_3|^4+|z_3-z_1|^4}{(|z_1{-}z_2|^2+|z_2{-}z_3|^2+|z_3{-}z_1|^2)^2}.$$
\end{lemma}

\begin{proof} The equality in the lemma is trivial if $u=0$. So we assume that $u>0$. For every $k\in\{1,2,3\}$ write the complex number $z_k$ as $x_k+iy_k$ for some real numbers $x_k,y_k$. Since $|\langle \vec z-z_0|\trig\rangle|$ is invariant under isometric transformations of the triangle $\vec z$, we lose no generality assuming that $z_1=0$ and $z_2=x_2$ is a positive real number. 
In this case $\sign(\vec z\,)=\sign(y_3)$.

Denote the lengths of the sides of the triangle $\vec z$ by $a=|z_1-z_2|$, $b=|z_2-z_3|$ and $c=|z_3-z_1|$. Let $\beta$ be the angle of the triangle $\vec z$ at the vertex $z_1=0$. So, $\beta$ is opposite to the side of length $b$. If $z_3=0$, then we put $\beta=0$. 
 
It follows that $x_3=c\cdot \cos(\beta)$ and $y_2=\sign(\vec z\,)\cdot c\cdot \sin(\beta)$. By the cosine theorem, $b^2=|z_2-z_3|^2=a^2+c^2-2ac\cos(\beta)$ and hence
$$\cos(\beta) = \frac{a^2+c^2-b^2}{2ac}$$
and
$$
\begin{aligned}
\sin(\beta) &=\sqrt{1-\Big( \frac{a^2+c^2-b^2}{2ac}\Big)^2} =\\
&=\frac{ \sqrt{2a^2b^2+2a^2c^2+2b^2c^2-(a^4+b^4+c^4)}}{2ac} = \frac{ \sqrt{(a^2+b^2+c^2)^2-2(a^4+b^4+c^4)}}{2ac}=\\
&=\frac{3u^2\sqrt{1-2q}}{2ac}
\end{aligned}
$$

By Lemma~\ref{l2}, 
$$	
   \begin{aligned}
   |\langle \vec z-z_0|\trig\rangle|^2&=\tfrac13\big(x_2^2+ x_3^2 + y_3^2 - x_2x_3+ \sqrt{3}x_2y_3\big) =
    \tfrac13\big(a^2+c^2 -ac\cos (\beta)+\sign(\vec z\,) \sqrt{3}ac\sin(\beta)\big)= \\
    &=\tfrac16\big(a^2+b^2+c^2+\sign(\vec z\,)3u^2\sqrt{3-6q})=\\
&=\frac{u^2}2\big(1+\sign(\vec z\,)\sqrt{3-6q}\,\big).
  \end{aligned}
  $$
\end{proof}

Now we are able to prove the main result of this section.

\begin{theorem}\label{t:main} For any  triangle $\vec z=(z_1,z_2,z_3)\in\IC^3$ on the complex plane 
\begin{enumerate}
\item the isometric deviation $d(\IT\vec z+\IC,\trig)=\sqrt{1+u^2-\sqrt{2}u\sqrt{1+\sign(\vec z\,)\sqrt{3-6q}}}$;
\item the affine deviation $d(\IC^*\vec z+\IC,\trig)=\sqrt{\frac12\big(1-\sign(\vec z\,)\sqrt{3-6q}\,\big)}$,
\end{enumerate}
where $$u=\sqrt{\frac{|z_1-z_2|^2+|z_2-z_3|^2+|z_3-z_1|^2}3}\quad\mbox{and}\quad q=\frac{|z_1-z_2|^4+|z_2-z_3|^4+|z_3-z_1|^4}{\big(|z_1{-}z_2|^2+|z_2{-}z_3|^2+|z_3{-}z_1|^2\big)^2}.$$
\end{theorem}

\begin{proof} Let $z_0=\frac13(z_1+z_2+z_3)$. By Theorem~\ref{t5} and Lemmas~\ref{l1}, \ref{l3}, we have  
$$
d(\IT\vec z+\IC,\trig)=
\sqrt{\|\vec z-z_0\|^2+1-2\cdot|\langle \vec z-z_0|\trig\rangle|}=\sqrt{1+u^2-\sqrt{2}u\sqrt{1+\sign(\vec z\,)\sqrt{3-6q}}}
$$
and
$$
d(\IC^*\vec z+\IC,\trig)=\sqrt{1-\frac{|\langle \vec z-c|\trig\rangle|^2}{\|\vec z-c\|^2}}=\sqrt{1-\frac{u^2(1+\sign(\vec z\,)\sqrt{3-6q}\,\big)}{2u^2}}=\sqrt{\tfrac12\big(1-\sign(\vec z\,)\sqrt{3-6q}\,\big)}.
$$
\end{proof}

\section{Interplay between the affine deviation and the unbalance factor}\label{s:K}

In this section we investigate the interplay between the unbalance factor $\unbal$ of a triangle $\vec z=(z_1,z_2,z_3)$ on the complex plane and the affine deviation $d(\IC^*\vec z+\IC,\trig)$ of the triangle $\vec z$ from the regular triangle $\trig$.

We recall that 
$$\unbal=\frac{|z_1+z_2e^{-i\frac{2\pi}3}+z_3 e^{i\frac{2\pi}3}|}{|z_1+z_2e^{i\frac{2\pi}3}+z_3 e^{-i\frac{2\pi}3}|}.$$
The unbalance factor $\unbal$ was introduced by Fortescue \cite{Fort} and is widely used in Electric Engineering \cite{Blac}, \cite{Das}, \cite{IEEE} for evaluation of the quality of 3-phase electric power.

\begin{theorem}\label{t:K} Any triangle $\vec z=(z_1,z_2,z_3)$ on the complex plane has unbalance factor
$$\unbal=\sqrt{\frac{1-\sign(\vec z\,)\sqrt{3-6q}}{1+\sign(\vec z\,)\sqrt{3-6q}}}$$where $q$ is the quadrofactor of the triangle $\vec z$.
\end{theorem}

\begin{proof} 
Let $z_0=\frac13(z_1+z_2+z_3)$ be the center of the triangle $\vec z$ and $s=|z_1-z_2|^2+|z_2-z_3|^2+|z_3-z_1|^2$. Observe that 
the positive component $\frac13(z_1+z_1e^{-i\frac{2\pi}3}+z_2e^{i\frac{2\pi}3})$ of $\vec z$ has length
$$
\begin{aligned}
\tfrac13|z_1+z_2e^{-i\frac{2\pi}3}+z_3e^{i\frac{2\pi}3}|&=\tfrac1{\sqrt{3}}|\langle \vec z\,|\trig\rangle|=\tfrac1{\sqrt{3}}|\langle \vec z-z_0|\trig\rangle+z_0(1+e^{-i\frac{2\pi}3}+e^{i\frac{2\pi}3})|=\\
&=\tfrac1{\sqrt{3}}|\langle \vec z-z_0|\trig\rangle+0|=\tfrac1{\sqrt{3}}|\langle \vec z-z_0|\trig\rangle|=\tfrac{1}{\sqrt{6}}u\sqrt{1+\sign(\vec z\,)\sqrt{3-6q}},
\end{aligned}
$$
according to Lemma~\ref{l3}.

For evaluation of the length of the negative component $z_1+z_1e^{i\frac{2\pi}3}+z_2e^{-i\frac{2\pi}3}$, consider the triangles
$$\cev z=(\bar z_1,\bar z_2,\bar z_3),\quad\bar\trig=\tfrac1{\sqrt{3}}(1,e^{-i\frac{2\pi}3},e^{i\frac{2\pi}3})$$ and observe that 
$$
\begin{aligned}
\tfrac13|z_1+z_2e^{i\frac{2\pi}3}+z_3e^{-i\frac{2\pi}3}|&=\tfrac1{\sqrt{3}}|\langle \vec z\,|\bar\trig\rangle|=\tfrac1{\sqrt{3}}|\langle \vec z-z_0|\bar \trig\rangle+ z_0(1+e^{i\frac{2\pi}3}+e^{-i\frac{2\pi}3})|=\\
&=\tfrac1{\sqrt{3}}|\langle \vec z-z_0|\bar\trig\rangle+0|=\tfrac13|\langle \vec z-z_0|\trig\rangle|=\tfrac1{\sqrt{3}}|\langle \cev z-\bar z_0|\trig\rangle|=\\
&=\tfrac1{\sqrt{6}}u\sqrt{1+\sign(\cev z\,)\sqrt{3-6q}}=\tfrac1{\sqrt{6}}u\sqrt{1-\sign(\vec z\,)\sqrt{3-6q}},
\end{aligned}
$$
according to Lemma~\ref{l3}.

Then 
$$\unbal=\frac{|z_1+z_2e^{i\frac{2\pi}3}+z_3e^{-i\frac{2\pi}3}|}{|z_1+z_2e^{-i\frac{2\pi}3}+z_3e^{i\frac{2\pi}3}|}=\frac{u\sqrt{\tfrac16\big(1-\sign(\vec z\,)\sqrt{3-6q}\,\big)}}{u\sqrt{\tfrac16\big(1+\sign(\vec z\,)\sqrt{3-6q}\,\big)}}=\frac{\sqrt{1-\sign(\vec z\,)\sqrt{3-6q}}}{\sqrt{1+\sign(\vec z\,)\sqrt{3-6q}}}.$$
\end{proof}

Theorems~\ref{t:K} and \ref{t:main} imply the following corollary expressing the unbalance factor of a triangle via its affine deviation for the regular triangle $\trig$ and vice versa.

\begin{corollary}\label{c:K} For any triangle $\vec z$ in the complex plane 
$$\unbal=\sqrt{\frac{d(\IC^*\vec z+\IC,\trig)^2}{1-d(\IC^*\vec z+\IC,\trig)^2}}\quad\mbox{and}\quad d(\IC^*\vec z+\IC,\trig)=\sqrt{\frac{\unbal^2}{1+\unbal^2}}.$$
\end{corollary}

\begin{proof} By Theorem~\ref{t:main},
$$
d(\IC^*\vec z+\IC,\trig)^2=\tfrac12\big(1-\sign(\vec z\,)\sqrt{3-6q}\,\big),$$ 
and hence $$\sign(\vec z\,)\sqrt{3-6q}=1-2d(\IC^*\vec z+\IC,\trig)^2.$$
After substitution of $\sign(\vec z\,)\sqrt{3-6q}$ into the formula for $\unbal$ from Theorem~\ref{t:K}, we obtain the desired equality
$$
\unbal=\sqrt{\frac{1-\sign(\vec z\,)\sqrt{3-6q}}{1+\sign(\vec z\,)\sqrt{3-6q}}}=\sqrt{\frac{d(\IC^*\vec z+\IC,\trig)^2}{1-d(\IC^*\vec z+\IC,\trig)^2}},
$$
which implies 
$$
\unbal^2\cdot\big(1-d(\IC^*\vec z+\IC,\trig)^2\big)=d(\IC^*\vec z+\IC,\trig)^2
$$
and finally
$$d(\IC^*\vec z+\IC,\trig)^2=\frac{\unbal^2}{1+\unbal^2}.$$
\end{proof}

\section{Visualization of the isometric and affine deviations}

For applications to problems of quality control in Electric Engineering, it is important to visualize the isometric and affine deviations of a given triangle with sides $a,b,c$ from the unit triangle $\triangle$. This can be done using the following theorem.

\begin{theorem} Let $\vec z=(z_1,z_2,z_3)\in\IC^3$ be a non-singular plane triangle with sides $a=|z_1-z_2|$, $b=|z_2-z_3|$ and $c=|z_3-z_1|$. 
Let $$u=\sqrt{\frac{a^2+b^2+c^2}3},\quad q=\frac{a^4+b^4+c^4}{(a^2+b^2+c^2)^2},$$
$$\begin{aligned}
z_1^\circlearrowleft&=\frac{2a^2+2c^2-b^2+3\sign(\vec z)u^2\sqrt{3-6q}}{6\sqrt{3}}+\frac{c^2-a^2}{6}\cdot i,\\
z_2^\circlearrowleft&=\frac{c^2-5a^2+b^2-3\sign(\vec z)u^2\sqrt{3-6q}}{12\sqrt{3}}+\frac{a^2-c^2+3b^2+3\sign(z)u^2\sqrt{3-6q}}{12}\cdot i,\\
z_3^\circlearrowleft&=\frac{a^2-5c^2+b^2-3\sign(\vec z)u^2\sqrt{3-6q}}{12\sqrt{3}}-\frac{c^2-a^2+3b^2+3\sign(\vec z)u^2\sqrt{3-6q}}{12}\cdot i.
\end{aligned}
$$
Then 
$$d(\IC^*\vec z+\IC,\triangle)=d(\vec z_*,\triangle)\quad\mbox{and}\quad d(\IT\vec z+\IC,\triangle)=d(\vec z_{\star},\triangle),$$
where 
$$\vec z_*=\frac{1}{u^2}\cdot(z_1^\circlearrowleft,z_2^\circlearrowleft,z_3^\circlearrowleft)\quad\mbox{and}\quad\vec z_\star=
\frac{\sqrt{2}}{u\sqrt{1+\sign(\vec z)\sqrt{3-6q}}}\cdot(z_1^\circlearrowleft,z_2^\circlearrowleft,z_3^\circlearrowleft).$$
If $a=b=c$ and $\sign(\vec z)=-1$, then $1+\sign(\vec z\,)\sqrt{3-6q}=0$ and the last formula does not determine $\vec z_\star$. In this case we can take $\vec z_\star=a\cdot\bar\trig$.
\end{theorem}

\begin{proof} For every $k\in\{1,2,3\}$, write the complex number $z_k$ as $x_k+iy_k$ for some real number $x_k,y_k$. Since the affine and isometric deviations of the triangle $\vec z$ from $\trig$ are invariant under isometric moves of the triangle $\vec z$, we can assume that $z_1=0$ and $z_2$ is a non-negative real number. So, $z_2=x_2\ge 0$. Let $\beta$ be the angle of the triangle $\vec z$ at the vertex $z_1$ (which is opposite to the side $|z_2-z_3=b$). It follows that $x_1=y_1=z_1=0$, $x_2=z_2=|z_2-z_1|=a$, $x_3=|z_1-z_3|\cos(\beta)=c\cos(\beta)$ and $y_3=\sign(\vec z)c\sin(\beta)$. It will be convenient to denote the number $\sign(\vec z)\in\{-1,0,1\}$ by $\pm$ and $-\sign(\vec z)$ by $\mp$. Then $y_3=\pm c\sin(\beta)$.

By the cosine theorem, $$\cos(\beta)=\frac{a^2+c^2-b^2}{2ac}.$$Then
$$
\begin{aligned}
\sin(\beta)&=\sqrt{1-\frac{(a^2+c^2-b^2)^2}{4a^2c^2}}=\frac{\sqrt{2(a^2c^2+b^2c^2+a^2b^2)-(a^4+b^4+c^4)}}{2ac}=\\
&=\frac{\sqrt{(a^2+b^2+c^2)^2-2(a^4+b^4+c^4)}}{2ac}=\frac{(a^2+b^2+c^2)\sqrt{1-2q}}{2ac}=\\
&=\frac{3u^2\sqrt{1-2q}}{2ac}=\frac{u^2\sqrt{3-6q}}{2ac}.
\end{aligned}
$$

 Let $$z_0=\tfrac13(z_1+z_2+z_3)=\tfrac{x_2+x_3}3+\tfrac{y_3}3i$$be the center of the triangle $\vec z$. Then 
$$\vec z-z_0=\big(-\tfrac{x_2+x_3}3-\tfrac{y_3}3i,\tfrac{2x_2-x_3}{3}-\tfrac{y_3}3i,\tfrac{2x_3-x_2}3+\tfrac{2y_3}3i\big).$$
In the proof of Lemma~\ref{l2} we have derived the formula:
$$\langle \vec z-z_0|\trig\rangle=-\tfrac{x_2+x_3+\sqrt 3 y_3}{2\sqrt{3}}-\tfrac{\sqrt{3}x_2-\sqrt 3x_3+y_3}{2\sqrt{3}}\cdot i.$$
Then 
$$
\begin{aligned}
&\overline{\langle z-z_0|\trig\rangle}\cdot(\vec z-z_0)=\\
&=\big(-\tfrac{x_2+x_3+\sqrt 3 y_3}{2\sqrt{3}}+\tfrac{\sqrt{3}x_2-\sqrt 3x_3+y_3}{2\sqrt{3}}\cdot i\big)\cdot\big(-\tfrac{x_2+x_3}3-\tfrac{y_3}3i,\tfrac{2x_2-x_3}{3}-\tfrac{y_3}3i,\tfrac{2x_3-x_2}3+\tfrac{2y_3}3i\big)=\\
&=\big(\tfrac{(x_2+x_3+\sqrt 3 y_3)(x_2+x_3)+(\sqrt{3}x_2-\sqrt 3x_3+y_3)y_3}{6\sqrt{3}}+\tfrac{(x_2+x_3+\sqrt 3 y_3)y_3-(\sqrt{3}x_2-\sqrt 3x_3+y_3)(x_2+x_3)}{6\sqrt{3}}i,\\
&\hskip20pt\tfrac{-(x_2+x_3+\sqrt 3 y_3)(2x_2-x_3)+\sqrt{3}x_2y_3-\sqrt 3x_3y_3+y_3^2}{6\sqrt{3}}+\tfrac{(\sqrt{3}x_2-\sqrt 3x_3+y_3)(2x_2-x_3)+(x_2+x_3+\sqrt 3 y_3)y_3}{6\sqrt{3}}i,\\
&\hskip20pt\tfrac{-(x_2+x_3+\sqrt 3 y_3)(2x_3-x_2)-2\sqrt{3}x_2y_3+2\sqrt 3x_3y_3-2y_3^2}{6\sqrt{3}}+\tfrac{(\sqrt{3}x_2-\sqrt 3x_3+y_3)(2x_3-x_2)-2(x_2+x_3+\sqrt 3 y_3)y_3}{6\sqrt{3}}i\big)=\\
&=\big(\tfrac{x_2^2+x_3^2+2x_2x_3+2\sqrt{3}x_2y_3+y_3^2}{6\sqrt{3}}+\tfrac{x_3^2-x_2^2+y_3^2}{6}i,\\
&\hskip20pt\tfrac{x_3^2-2x_2^2-x_2x_3-\sqrt 3 x_2y_3+y_3^2}{6\sqrt{3}}+\tfrac{2x_2^2+x_3^2-3x_2x_3+\sqrt{3}x_2y_3+y_3^2}{6}i,\\
&\hskip20pt\tfrac{x_2^2-2x_3^2-x_2x_3-\sqrt{3}x_2y_3-2y_3^2}{6\sqrt{3}}-\tfrac{x_2^2+2x_3^2-3x_2x_3+\sqrt{3}x_2y_3+2y_3^2}{6}i\big)=\\
&=\big(\tfrac{a^2+c^2\cos(\beta)^2+2ac\cos(\beta)\pm2\sqrt{3}ac\sin(\beta)+c^2\sin(\beta)^2}{6\sqrt{3}}+\tfrac{c^2\cos(\beta)^2-a^2+c^2\sin(\beta)^2}{6}i,\\
&\hskip20pt\tfrac{c^2\cos(\beta)^2-2a^2-ac\cos(\beta)\mp\sqrt 3 ac\sin(\beta)+c^2\sin(\beta)^2}{6\sqrt{3}}+\tfrac{2a^2+c^2\cos(\beta)^2-3ac\cos(\beta)\pm\sqrt{3}ac\sin(\beta)+c^2\sin(\beta)^2}{6}i,\\
&\hskip20pt\tfrac{a^2-2c^2\cos(\beta)^2-ac\cos(\beta)\mp\sqrt{3}ac\sin(\beta)-2c^2\sin(\beta)^2}{6\sqrt{3}}-\tfrac{a^2+2c^2\cos(\beta)^2-3ac\cos(\beta)\pm\sqrt{3}ac\sin(\beta)+2c^2\sin(\beta)^2}{6}i\big)=\\
&=\big(\tfrac{a^2+c^2+2ac\cos(\beta)\pm2\sqrt{3}ac\sin(\beta)}{6\sqrt{3}}+\tfrac{c^2-a^2}{6}i,\\
&\hskip20pt\tfrac{c^2-2a^2-ac\cos(\beta)\mp\sqrt 3 ac\sin(\beta)}{6\sqrt{3}}+\tfrac{2a^2+c^2-3ac\cos(\beta)\pm\sqrt{3}ac\sin(\beta)}{6}i,\\
&\hskip20pt\tfrac{a^2-2c^2-ac\cos(\beta)\mp\sqrt{3}ac\sin(\beta)}{6\sqrt{3}}-\tfrac{a^2+2c^2-3ac\cos(\beta)\pm\sqrt{3}ac\sin(\beta)}{6}i\big)=\\
&=\big(\tfrac{a^2+c^2+(a^2+c^2-b^2)\pm3\sqrt{3}u^2\sqrt{1-2q})}{6\sqrt{3}}+\tfrac{c^2-a^2}{6}i,\\
&\hskip20pt\tfrac{c^2-2a^2-\tfrac12(a^2+c^2-b^2)\mp\tfrac32\sqrt 3u^2\sqrt{1-2q}}{6\sqrt{3}}+\tfrac{2a^2+c^2-\tfrac32(a^2+c^2-b^2)\pm\tfrac32\sqrt{3}u^2\sqrt{1-2q}}{6}i,\\
&\hskip20pt\tfrac{a^2-2c^2-\tfrac12(a^2+c^2-b^2)\mp\tfrac32\sqrt{3}u^2\sqrt{1-2q}}{6\sqrt{3}}-\tfrac{a^2+2c^2-\tfrac32(a^2+c^2-b^2)\pm\tfrac32\sqrt{3}u^2\sqrt{1-2q}}{6}i\big)=\\
&=\big(\tfrac{2a^2+2c^2-b^2\pm3u^2\sqrt{3-6q})}{6\sqrt{3}}+\tfrac{c^2-a^2}{6}i,\\
&\hskip20pt\tfrac{c^2-5a^2+b^2\mp 3u^2\sqrt{3-6q}}{12\sqrt{3}}+\tfrac{a^2-c^2+3b^2\pm3u^2\sqrt{3-6q}}{12}i,\\
&\hskip20pt\tfrac{a^2-5c^2+b^2\mp3u^2\sqrt{3-6q}}{12\sqrt{3}}-\tfrac{c^2-a^2+3b^2\pm3u^2\sqrt{3-6q}}{12}i\big)=\\
&=(z_1^\circlearrowleft,z_2^\circlearrowleft,z_3^\circlearrowleft).
\end{aligned}
$$
Applying Theorem~\ref{t5} and Lemmas~\ref{l1}, \ref{l3}, we conclude that
$$d(\IC^*\vec z+\IC,\trig)=d(\vec z_*,\trig)\quad\mbox{and}\quad d(\IT\vec z+\IC,\trig)=d(\vec z_\star,\trig),$$
where
$$\vec z_*=\frac{\overline{\langle\vec z-z_0|\trig\rangle}(\vec z-z_0)}{\|\vec z-z_0\|^2}=\frac{(z_1^\circlearrowleft,z_2^\circlearrowleft,z_3^\circlearrowleft)}{u^2}
$$
and
$$
\vec z_\star=\frac{\overline{\langle\vec z-z_0|\trig\rangle}(\vec z-z_0)}{|\langle \vec z-z_0|\trig\rangle|}=\frac{(z_1^\circlearrowleft,z_2^\circlearrowleft,z_3^\circlearrowleft)}{\tfrac1{\sqrt{2}}u\sqrt{1+\sign(\vec z)\sqrt{3-6q}}}.
$$
If $a=b=c$ and $\sign(\vec z)=-1$, then $1+\sign(\vec z\,)\sqrt{3-6q}=0$ and the last formula does not determine $z_\star$. In this case Theorem~\ref{t2} allows us to take for $\vec z_\star$ any triangle $t(\vec z-z_0)$ where $t\in\IT$. So, we can choose $t\in\IT$ such that $\vec z_\star=t(\vec z-z_0)=a\cdot\bar \trig$.
\end{proof}

\section{Acknowledgements}

The authors express their sincere thanks to the {\tt Mathoverflow} user {\tt \@bathhalf15320} for suggesting a formula connecting the quadrofactor and  normalized area of a plane triangle in Theorem~\ref{t:Heron}, see {\tt https://mathoverflow.net/a/389681/61536}.

\end{document}